\documentstyle[amstex,amscd]{article}
\title{A note on symplectic singularities}
\author{Yoshinori Namikawa}
\date{ }
\begin{document}
\maketitle

\begin{center}
{\bf Introduction}
\end{center}

Let $X$ be an affine (or Stein) open subset 
of an algebraic variety 
of $\dim n$ over ${\bold C}$. $X$ is called a 
symplectic singularity if the smooth locus $U$ 
admits an everywhere non-degenerate 
holomorphic closed 2-form $\omega$ and, 
for an arbitrary resolution 
$f: Y \to X$ with $f^{-1}(U) \cong U$, 
$\omega$ extends to a regular 2-form on $Y$. Let 
$\Sigma$ be the singular locus. 
In this note we shall prove the following. 
\vspace{0.2cm}

{\bf Theorem}. {\em Let $X$ be a 
symplectic singualrity. Then $\Sigma$ has 
no codimension 3 irreducible components.} 
\vspace{0.2cm}

{\bf Corollary 1}. {\em Let $X$ be a symplectic 
singularity. Then $X$ is terminal if and only 
if $\mathrm{Codim}(\Sigma \subset X) \geq 4$ 
}. \vspace{0.15cm}

{\em Proof of Corollary 1}. Assume that $X$ has 
only terminal singularities. Then  
$\mathrm{Codim}(\Sigma \subset X) \geq 3$. 
By Theorem, $\Sigma$ has no codimension 3 
irreducible components. 
Therefore 
$\mathrm{Codim}(\Sigma \subset X) \geq 4$.

Conversely, assume 
$\mathrm{Codim}(\Sigma \subset X) \geq 4$. 
Take a resolution $f: Y \to X$ of $X$ and 
denote by $\{E_i\}$ the $f$-exceptional 
divisors. Since $X$ has canonical singularities, 
we have $K_Y = f^*K_X + \Sigma a_iE_i$ with 
non-negative integers $a_i$. One can prove that 
$a_i$ are actually all positive by the same 
argument as \S 1 \footnote{Although it is 
not explicitly mentioned, the following fact is 
proved in \S 1: Let $X$ be a symplectic singularity 
and $f: Y \to X$ a resolution of 
singularities. Assume that $E_i$ is an irreducible 
$f$-exceptional divisor such that $a_i = 0$. 
Then $\dim f(E_i) = n-2$, where $n = \dim X$.}.  

\vspace{0.12cm}

We can define a symplectic variety $Z$ 
similarly as a compact normal Kaehler space $Z$ 
whose regular locus $V$ admits an everywhere 
non-degenerate holomorphic closed 2-form 
which extends to a regular 2-form on an 
arbitrary resolution $\tilde Z$ of $Z$. 
Then the following result is also 
a corollary of our theorem.
\vspace {0.15cm}

{\bf Corollary 2}. {\em Let $Z$ be a projective 
symplectic variety with terminal singularities. 
Then the Kuranishi space $\mathrm{Def}(Z)$ 
is smooth.} \vspace{0.15cm}

{\em Proof of Corollary 2}. We already know 
the result when $Z$ is a symplectic variety 
with $\mathrm{Codim}(\Sigma \subset Z) \geq 4$ 
by [Na 2, Theorem (2.4)]. Now the 
result follows from Corollary 1.
Q.E.D. \vspace{0.2cm} 

In the remainder we put $\dim X = n$. 
Our strategy is the following. 
Assume that $\Sigma$ has an $n-3$ 
dimensional irreducible component $\Sigma'$. 
Let $H : = H_1 \cap H_2 \cap ... \cap H_{n-3}$ 
be a complete intersection of $(n-3)$ 
general hyperplane sections of $X$. 
$H$ intersects $\Sigma'$ in a 
finite number of points. 
Pick a point $p$ from these points and 
we shall only consider 
a small open neighborhood of $p \in X$. 
For a resolution $f: Y \to X$, 
$\tilde H := f^{-1}(H)$ is a 
resolution of $H$ by the theorem of Bertini. 
Since $X$ has canonical singularities, 
we have $K_Y = f^*K_X + \Sigma a_iE_i$ 
where $E_i$ are $f$-exceptional divisors 
and the coefficients $a_i$ are all 
{\em non-negative}. By the adjunction formula, 
$K_{\tilde H} = 
f^*K_H + \Sigma a_i(E_i \cap \tilde H)$. 
As $F_i := E_i \cap \tilde H$ are exceptional 
divisors of $f\vert_{\tilde H}$, 
we know that $H$ has also canonical singularities.  
By taking an (n-3)-parameter family 
of complete intersections $H_s$ with 
$H_0 = H$, we can regard $X$, locally around $p$, 
as the total space ${\cal H}$ of the 
(n-3)-parameter family. 
By definition, there is a map 
$\pi: {\cal H} \to \Delta^{n-3}$ 
such that each fiber is $H_s$. 

If the 3-dimensional canonical singularity 
$(H, p)$ is not compound Du Val (cDV) singularity, 
then there is at least one 
$f\vert_{\tilde H}$-exceptional divisor  
$F_i$ over $p$ such that $a_i = 0$ by 
[Re, Theorem 2.2]. 
Therefore, there is at least one 
$f$-exceptional divisor $E_i$ with $a_i = 0$. 
This $E_i$ is mapped onto $\Sigma'$ by $f$. 
We shall show that this contradicts the fact 
that $X$ is a symplectic singularity in \S 1.

Next consider the case where $(H, p)$ 
is a cDV point. $({\cal H}, p)$ is a 
hypersurface singularity in $(\Delta^{n+1}, 0)$ 
defined by $g(x,y,z,t,s_1,s_2, ..., s_{n-3}) = 0$, 
where $(x, y, z, t, s_1, s_2, ..., s_{n-3})$ 
are coordinates of $\Delta^{n+1}$, 
and ${\cal H} \to \Delta^{n-3}$ is given by 
$(x, y, z, t, s_1, s_2, ..., s_{n-3}) 
\to (s_1, s_2, ..., s_{n-3})$. 
Since $({\cal H}, p)$ is an (n-3)-parameter 
deformation of an isolated cDV point $(H, p)$ 
and since $(H, p)$ is a one parameter 
deformation of a Du Val singularity $S$, 
$({\cal H}, p)$ is an (n-2)-parameter deformation 
of a Du Val singularity $S$. 
We may assume that the n-2 parameters are given by 
$t$ and $s_1, s_2, ..., s_{n-3}$. 
Let $\mathrm{Def}(S)$ be the semi-universal 
deformation space (=Kuranishi space) of $S$. 
Then we have a holomorphic map 
$\phi: \Delta^{n-2} \to \mathrm{Def}(S)$ 
such that $\phi(0) = 0 \in \mathrm{Def}(S)$. 
Let $\omega$ be a holomorphic 2-form defined 
on the regular part of ${\cal H}$. 
By using the above description of ${\cal H}$, 
we shall write $\omega$ rather explicitly 
in terms of $dx$, $dy$, $dz$, $dt$, 
$ds_1$, ..., $ds_{n-3}$. 
This explicit description of $\omega$ 
will tell us that $\wedge^{n/2}\omega$ 
is somewhere-vanishing (cf. \S 2).  
\vspace{0.12cm}

{\bf Remark}. Since a symplectic singularity 
is a Gorenstein 
canonical singularity, 
there is a closed subset $\Sigma_0$ of 
$\Sigma$ such that each point 
$p \in \Sigma \setminus \Sigma_0$ 
has an open neighborhood (in $X$) 
isomorphic to $(R.D.P) \times 
({\bold C}^{n-2}, 0)$ ([Re]). 
Then our theorem can be generalized 
to the following form: \vspace{0.12cm}

{\em  Let $X$ be a symplectic singularity. 
Then $\mathrm{Codim}(\Sigma_0 \subset X) \geq 4$.} 
\vspace{0.12cm}

The proof is almost the same as Theorem. 
When $X$ has a 
symplectic resolution, 
this result is proved in [Na 2] 
by a different method. \vspace{0.2cm}  

\begin{center}
{\bf \S 1. non-cDV case}
\end{center}

We use the same notation as Introduction. 
This section deals with the 
case $(H,p)$ is not a cDV point. 
As in Introduction, for a resolution 
$f: Y \to X$, there is an 
$f$-exceptional divisor $E_i$ with $a_i = 0$ 
which dominates $\Sigma'$. 
We put ${\tilde{\cal H}} := f^{-1}({\cal H})$. 
Denote by the same notation $f$ (resp. $E_i$) 
the restriction of $f$ 
to ${\tilde{\cal H}}$ 
(resp. the restriction of 
$E_i$ to ${\tilde{\cal H}}
$). We shall derive a contradiction 
by assuming that there is an 
everywhere non-degenerate 2-form $\omega$ 
on the smooth locus of $\cal 
H$.

Since $a_i = 0$, $\wedge^{n/2}f^*\omega$ 
does not vanish at general 
point $q \in E_i$. This means that 
$f^*\omega$ is a non-degenerate 
2-form at $q$. We can choose a system 
of local coordinates $(x_1, x_2, 
..., x_n)$ at $q \in \tilde{\cal H}$ 
in such a way that $x_1 = 0$ is a 
defining equation of $E_i$ and 
$f^*\omega(q) \in H^0(\tilde{\cal H}, 
\Omega^2_{\tilde{\cal H}}\otimes k(q))$ 
is of the following form 

$$  dx_1 \wedge dx_2 + dx_3 
\wedge dx_4 + ... + dx_{n-1}\wedge dx_n. 
$$ 

Let $\omega' \in H^0(E_i, \Omega^2_{E_i})$ 
be the restriction of 
$f^*\omega$ to $E_i$. 
We see by the above expression of $f^*\omega(q)$  
that $\wedge^{n/2-1}\omega' \ne 0$ at $q$.  

Take a birational proper morphism 
$\nu: {\cal W} \to \tilde{\cal H}$ 
from a smooth variety ${\cal W}$ to 
$\tilde{\cal H}$ in such a way that 
$F := (f\circ \nu)^{-1}(\Sigma')$ is a 
simple normal crossing divisor 
of ${\cal W}$. Denote by $\hat\Omega^i_F$ 
(resp.$\hat\Omega^i_{F/\Sigma'}$)  
the sheaf of i-forms of $F$ 
modulo torsion (resp. the sheaf of 
relative i-forms of $F$ over 
$\Sigma'$ modulo torsion). 
Let $\omega'' \in H^0(F, \hat\Omega^2_F)$ 
be the restriction of 
$(f\circ\nu)^*\omega$ to $F$. We may assume that 
$\nu$ is an isomorphism at $q$, and 
$\omega'' = \omega'$ at $q$. 
Therefore, $\wedge^{n/2-1}\omega'' \ne 0$ at $q$.  
We conclude that 
$\omega''$ is not the pull-back of 
any 2-form of $\Sigma'$ by the map 
$F \to \Sigma'$ because 
$\dim \Sigma' \leq n-3$.
 Denote by $F_x$ the 
fiber of the map $F \to \Sigma'$ 
over $x \in \Sigma'$.

We consider two exact sequences 

$$ 0 \to {\cal F} \to \hat\Omega^2_F 
\to \hat\Omega^2_{F/\Sigma'} \to 0, 
$$ 

$$ 0 \to (f\circ\nu)^*\Omega^2_{\Sigma'} 
\to {\cal F} \to 
(f\circ\nu)^*\Omega^1_{\Sigma'}
\otimes\hat\Omega^1_{F/\Sigma'} \to 
0.$$

We know that $H^0(F_x, \hat\Omega^2_{F_x}) 
\ne 0$ or $H^0(F_x, 
\hat\Omega^1_{F_x}) \ne 0$ for a general point 
$x \in \Sigma'$ by 
taking global sections of these 
exact sequences because $\omega''$ is 
not the pull-back of any 2-form of 
$\Sigma'$. By the mixed Hodge 
structure on $H^1(F_x)$ or $H^2(F_x)$, 
$H^1(F_x, {\cal O}_{F_x}) \ne 0$ 
if $H^0(F_x, \hat\Omega^1_{F_x}) \ne 0$ 
and $H^2(F_x, {\cal O}_{F_x}) 
\ne 0$ if 
$H^0(F_x, \hat\Omega^2_{F_x}) \ne 0$. 

On the other hand, let ${\cal H}_x$ be a 
general complete intersection 
of n-3 general hyperplanes passing 
through $x$. Put ${\cal W}_x := 
(f\circ\nu)^{-1}({\cal H}_x)$. 
Then the map 
$(f\circ\nu)\vert_{{\cal W}_x}: 
{\cal W}_x \to {\cal H}_x$ 
is a resolution of singularities, and 
$F_x$ is the 
exceptional locus of this resolution. 
Since ${\cal H}_x$ has only 
rational singularities, 
$R^i((f\circ\nu)\vert_{{\cal W}_x})_*{\cal 
O}_{{\cal W}_x} = 0$ for all $i > 0$. 
We deduce from this fact that 
$H^i(F_x, {\cal O}_{F_x}) = 0$ for all 
$i > 0$. This is a contradiction. \vspace{0.2cm}

\begin{center}
{\bf \S 2. cDV case}
\end{center}

We use the same notation as Introduction. 
This section deals with the case $(H, p)$ 
is a cDV point. 
Let $S$ be a general hyperplane section of 
$H$ passing through $p$. 
If $S$ is of type $A_n$ (resp. $D_n$, $E_n$), 
then we call $(H, p)$ is of type 
$cA_n$ (resp. $cD_n$, $cE_n$). 
As in Introduction, $({\cal H}, p)$ is a 
hypersurface singularity defined by 
$g(x,y,z,t,s_1, ..., s_{n-3}) = 0$ in 
$(\Delta^{n+1}, 0)$, where 
$(x, y, z, t, s_1, s_2, ..., s_{n-3})$ 
are coordinates of $\Delta^{n+1}$. 
The projection $\Delta^{n+1} \to \Delta^{n-2}$ 
($(x,y,z,t,s_1,s_2, ..., s_{n-3}) 
\to (t,s_1,s_2, ..., s_{n-3})$ 
induces a Du Val fibration 
${\cal H} \to \Delta^{n-2}$. 
$S$ is the central fiber of this fibration. 
Moreover, the projection 
$\Delta^{n+1} \to \Delta^{n-3}$ 
($(x,y,z,t,s_1,s_2, ..., s_{n-3}) 
\to (s_1,s_2, ..., s_{n-3})$) induces a 
cDV fibration ${\cal H} \to \Delta^{n-3}$. 
The central fiber of this fibration is $H$. 
We may assume that each fiber $H_s$ 
of this cDV fibration has a unique 
singular point at 
$(0, 0, 0, 0, s_1, ... s_{n-3})$ 
and its type as a cDV point is independent 
of $(s_1, ..., s_{n-3})$ (cf. Introduction). 
The Du Val fibration ${\cal H} \to \Delta^{n-2}$ 
induces a holomorphic map 
$\phi: \Delta^{n-2} \to \mathrm{Def}(S)$ 
such that $\phi(0) = 0 \in \mathrm{Def}(S)$. 
${\cal H}$ is isomorphic to the family of 
Du Val singularities over 
$\Delta^{n-2}$ obtained by pulling back 
the semi-universal family over 
$\mathrm{Def}(S)$ by $\phi$. 
Therefore $g(x,y,z,t,s_1, s_2, ..., s_{n-3})$ 
can be written in the following way 
according to the type of the 
Du Val singularity $S$ (cf. [Lo]): \vspace{0.25cm}

$A_n$: $ x^2 + y^2 +  z^{n+1} + 
\phi_1(t,s)z^{n-1} + ... + 
\phi_{n-1}(t,s)z + \phi_n(t,s)$ \vspace{0.12cm} 

$D_n (n \geq 4)$: $ x^2 + y^2z + z^{n-1} 
+ \phi_1(t,s)z^{n-2} + ... + 
\phi_{n-1}(t,s) + \phi_n(t,s)y$ \vspace{0.12cm}

$E_6$: $ x^2 + y^3 + z^4 + \phi_1(t,s)z^2 
+ ... + \phi_3(t,s) + y(\phi_4(t,s)z^2 
+ ... + \phi_6(t,s))$ \vspace{0.12cm}

$E_7$: $ x^2 + y^3 + yz^3 + 
y(\phi_1(t,s)z + \phi_2(t,s)) + 
\phi_3(t,s)z^4 + ... + \phi_7(t,s)$ \vspace{0.12cm}

$E_8$: $ x^2 + y^3 + z^5 + 
y(\phi_1(t,s)z^3 + ... + \phi_4(t,s)) + 
\phi_5(t,s)z^3 + ... + \phi_8(t,s)$ \vspace{0.25cm}

We abbreviated $s$ for $s_1, ..., s_{n-3}$. 
The $\phi_i(t,s)$'s satisfy the 
following properties, 
which will be used below: \vspace{0.12cm}

(a)  $\phi_i(0,s) \equiv 0$;

(b)  $\phi_i(t,s) \in 
t^2{\bold C} \{t, s_1, ..., s_{n-3}\}$ 
except the following cases \vspace{0.1cm}

$\phi_1(t,s)$ in the case $(D_n)$; 

$\phi_4(t,s)$ in the case $(E_6)$; 

$\phi_3(t,s)$ in the case $(E_7)$;

$\phi_1(t,s)$ in the case $(E_8)$.   
\vspace{0.15cm}

(a) is satisfied because the fibers of 
the Du Val fibration 
${\cal H} \to \Delta^{n-2}$ over 
$(0,s_1, ..., s_{n-3}) \in \Delta^{n-2}$ 
are Du Val singularities of the same type.  
We can check (b) in the following way. 
Assume that some other $\phi_i(t,s)$'s 
are not contained in 
$t^2{\bold C} \{t, s_1, ..., s_{n-3}\}$. 
Let $H_s$ be the fiber of the cDV fibration 
${\cal H} \to \Delta^{n-3}$ over 
$s = (s_1, ..., s_{n-3})$. Take $s$ general. 
Then the type of a cDV point $H_s$ 
does not coincide with the type of 
a Du Val singularity $S$. 
This contradicts our assumption. \vspace{0.15cm}

Let ${\cal H}^0$ and $S^0$ be the 
smooth locus of ${\cal H}$ and $S$ respectively. 
Denote by $p$ the Du Val fibration 
${\cal H} \to \Delta^{n-2}$. 
We write $p^0$ for the restriction of 
$p$ to ${\cal H}^0$. We have two exact sequences 

$$   0 \to {\cal F} \to 
\Omega^2_{{\cal H}^0} \to 
\Omega^2_{{\cal H}^0/\Delta^{n-2}} \to 0, $$

$$   0 \to (p^0)^*{\Omega}^2_{\Delta^{n-2}} \to 
{\cal F} \to 
\Omega^1_{{\cal H}^0/\Delta^{n-2}}
\otimes (p^0)^*\Omega^1_{\Delta^{n-2}} \to 0.  $$

Let $\omega$ be a holomorphic 2-form 
defined on ${\cal H}^0$, that is, 
$\omega \in H^0({\cal H}^0, \Omega^2_{{\cal H}^0})$. 
Let $\omega\vert_{S^0} \in 
H^0(S^0, \Omega^2_{{\cal H}^0}\vert_{S^0})$ 
be the image of $\omega$ by the natural map 
$H^0({\cal H}^0, \Omega^2_{{\cal H}^0}) 
\to H^0(S^0, \Omega^2_{{\cal H}^0}\vert_{S^0})$. 
We shall write $\omega\vert_{S^0}$ 
explicitly and compute 
$\wedge^{n/2}\omega\vert_{S^0}$.  

By restricting above exact sequences 
to $S^0$ and taking global sections, 
we have the exact sequences 

$$  0 \to H^0(S^0, {\cal F}\vert_{S^0}) 
\to H^0(S^0, \Omega^2_{{\cal H}^0}\vert_{S^0}) 
\stackrel{\alpha}\to H^0(S^0, \Omega^2_{S^0}), $$

$$  0 \to H^0(S^0, 
(p^0)^*\Omega^2_{\Delta^{n-2}}\vert_{S^0}) 
\to H^0(S^0, F\vert_{S^0}) \to 
H^0(S^0, \Omega^1_{S^0}\otimes 
(p^0)^*\Omega^1_{\Delta^{n-2}}). $$
As the first step we shall construct 
an explicit lift $\omega_1 
\in H^0(S^0, 
\Omega^2_{{\cal H}^0}\vert_{S^0})$ of 
$\alpha(\omega\vert_{S^0}) 
\in H^0(S^0, \Omega^2_{S^0})$.  \vspace{0.12cm}

{\bf Step 1}.  Consider the sheaf 
$\Omega^2_{{\cal H}^0/\Delta^{n-3}}$ of 
relative 2-forms with respect to the 
cDV fibration ${\cal H} \to \Delta^{n-3}$. 
By a generator $dy \wedge dz \wedge 
dt/(\partial g/ \partial x)$ of the relative 
dualizing sheaf $\omega_{{\cal H}/\Delta^{n-3}}$, 
the relative tangent sheaf 
$\Theta_{{\cal H}^0/\Delta^{n-3}}$ and 
$\Omega^2_{{\cal H}^0/\Delta^{n-3}}$ 
are identified. Thus we have 
$H^0({\cal H}^0, 
\Theta_{{\cal H}^0/\Delta^{n-3}}) 
\cong H^0({\cal H}^0, 
\Omega^2_{{\cal H}^0/\Delta^{n-3}})$.  
An element of $H^0({\cal H}^0, 
\Theta_{{\cal H}^0/\Delta^{n-3}})$ can be expressed as 

$$   f_1\partial /\partial x + 
f_2\partial/\partial y + 
f_3\partial/\partial z + f_4\partial/\partial t  $$  

with holomorphic functions $f_i$ on ${\cal H}$ 
which satisfy 
$f_1\partial g /\partial x + f_2\partial g/\partial y + 
f_3\partial g/\partial z + f_4\partial g/\partial t = 0$ 
on ${\cal H}$.  Therefore an element of 
$H^0({\cal H}^0, \Omega^2_{{\cal H}^0/\Delta^{n-3}})$ 
can be expressed as 

$$   f_2dz \wedge dt/(\partial g/\partial x) - 
f_3dy \wedge dt/(\partial g/\partial x) + 
f_4dy \wedge dz/(\partial g/\partial x).  $$

{\bf Lemma (2.1)}. {\em $f_i$ vanish at 
the origin for $i = 1, 2, 3, 4$.} \vspace{0.15cm}

{\em Proof}. We restrict the equation 
$f_1\partial g/\partial 
x + f_2\partial g/\partial y + 
f_3\partial g/\partial z + f_4\partial 
g/\partial t = 0$ on $\cal H$ to 
$H := \{s_1 = ... = s_{n-3} = 0\}$. 
 Note that 
$\partial g(x,y,z,t,s)/\partial x\vert_{s_1 = ... = 
s_{n-3} = 0} = \partial g(x,y,z,t,0)/\partial x$, 
and similar facts hold 
for the partial derivative of $y$, $z$ and $t$.  
$H$ is an isolated hypersurface singularity of 
${\bold C}^4 \ni (x,y,z,t)$ defined by 
$g(x,y,z,t,0) = 0$. Assume that some $f_i$ is 
unit. Then, applying Theorem of Frobenius, 
we can write on  ${\bold C}^4$ 

$$ f_1\partial /\partial x + 
f_2\partial /\partial y + f_3\partial 
/\partial z + f_4\partial /\partial t = 
(unit)\partial/\partial x' $$ 

by a suitable coordinates change of 
$(x,y,z,t)$ to $(x',y',z',t')$. 
Therefore, in ${\bold C}\{x',y',z',t'\}$

$$  \partial g/\partial x' + hg = 0  $$

with some $h \in {\bold C}\{x',y',z',t'\}$. 
Put $g' := g\cdot exp(\int 
hdx')$. Then we have 
$\partial g'/\partial x' = 0$. Since $g' = 0$ is 
a defining equation of $H$ in 
${\bold C}^4 \ni (x',y',z',t')$, this 
means that $H$ has non-isolated 
singularities or $H$ is smooth, which 
is a contradiction. Q.E.D. \vspace{0.2cm}
 
In particular, the image of 
$\omega$ by the map 
$H^0({\cal H}^0, \Omega^2_{{\cal H}^0}) 
\to H^0({\cal H}^0, 
\Omega^2_{{\cal H}^0/\Delta^{n-3}})$ 
has such an expression. 
We take $f_1, ..., f_4$ in such a way 
that the above 2-form is this image. 
Since $\Omega^2_{{\cal H}^0} \to 
\Omega^2_{{\cal H}^0/\Delta^{n-2}}$ 
factors through 
$\Omega^2_{{\cal H}^0/\Delta^{n-3}}$, 
$f_4dy \wedge dz/(\partial g/\partial x) 
\in H^0({\cal H}^0, 
\Omega^2_{{\cal H}^0/\Delta^{n-2}})$ 
coincides with the image of $\omega$ 
by the map $H^0({\cal H}^0, 
\Omega^2_{{\cal H}^0}) \to H^0({\cal H}^0, 
\Omega^2_{{\cal H}^0/\Delta^{n-2}})$.  
Therefore 
$\alpha(\omega\vert_{S^0}) = 
\overline{f}_4dy \wedge dz/(\partial g/\partial x) 
\in H^0(S^0, \Omega^2_{S^0})$, 
where $\overline{f}_4$ is the restriction of 
$f_4$ to $S^0$. 

Now the 2-form $\overline{f}_2dz 
\wedge dt/(\partial g/\partial x) - \overline{f}_3dy 
\wedge dt/(\partial g/\partial x) + \overline{f}_4dy 
\wedge dz/(\partial g/\partial x)$ 
can be considered as an element of 
$H^0(S^0, \Omega^2_{{\cal H}^0}\vert_{S^0})$ 
because $(\partial g/\partial s_i)\vert_{t=0} 
\equiv 0$ by the property (a) 
above\footnote{This is the reason why 
we consider $\omega\vert_{S^0}$ instead of 
$\omega$. In general, 
$f_2dz \wedge dt/(\partial g/\partial x) - 
f_3dy \wedge dt/(\partial g/\partial x) + 
f_4dy \wedge dz/(\partial g/\partial x)$ 
is not an element of 
$H^0({\cal H}^0, \Omega^2_{{\cal H}^0})$.}.  
So we put 

$$   \omega_1 := \overline{f}_2dz 
\wedge dt/(\partial g/\partial x) - 
\overline{f}_3dy \wedge 
dt/(\partial g/\partial x) + 
\overline{f}_4dy 
\wedge dz/(\partial g/\partial x). $$ 

{\bf Step 2}. Put 
$\omega_2 := \omega\vert_{S^0} - \omega_1$. 
By the construction of $\omega_1$, 
$\omega_2 \in H^0(S^0, {\cal F}\vert_{S^0})$. 
We shall write $\omega_2$ explicitly in 
terms of $dx$, $dy$, $dz$, $dt$, 
$ds_1$, ..., $ds_{n-3}$. 

By a generator 
$dy \wedge dz/(\partial g/\partial x)$ 
of the dualizing sheaf $\omega_S$, 
$\Theta_{S^0}$ and $\Omega^1_{S^0}$ 
are identified. Therefore 
$H^0(S^0, \Theta_{S^0}) \cong 
H^0(S^0, \Omega^1_{S^0})$. 
An element of $H^0(S^0, 
\Theta_{S^0})$ can be expressed as 

$$  h_1\partial/\partial x + 
h_2\partial/\partial y + 
h_3\partial/\partial z $$ 
with holomorphic functions $h_i$ on 
$S$ which satisfy 
$h_1\partial g/\partial x + 
h_2\partial g/\partial y + 
h_3\partial g/\partial z = 0$ on $S$. 
Therefore, an element of 
$H^0(S^0, \Omega^1_{S^0})$ 
can be expressed as 

$$  h_2 dz/(\partial g/\partial x) - 
h_3 dy/(\partial g/\partial x). $$ 

{\bf Lemma (2.2)} {\em $h_i$ vanish at 
the origin for $i = 1, 2, 3$}. \vspace{0.15cm}

{\em Proof}.  If some $h_i$ is a unit, 
then $S$ should have 
non-isolated singularities or 
$S$ should be smooth, by the same 
argument as the proof of Lemma (2.1). 
Q.E.D. \vspace{0.15cm}
 
The image of $\omega_2$ by the map 
$H^0(S^0, F\vert_{S^0}) \to 
H^0(S^0, \Omega^1_{S^0}\otimes 
(p^0)^*\Omega^1_{\Delta^{n-2}})$ has 
the expression

$$  \Sigma_{1 \leq j \leq n-3}\{h^{(j)}_2 
dz/(\partial g/\partial x) - 
h^{(j)}_3 dy/(\partial g/\partial x)\} 
\wedge ds_j  +  
\{h'_2 dz/(\partial g/\partial x) - 
h'_3 dy/(\partial g/\partial x)\} \wedge dt $$

for some $h^{(j)}_i$ and $h'_i$ 
which satify $h^{(j)}_1\partial g/\partial x + 
h^{(j)}_2\partial g/\partial y + 
h^{(j)}_3\partial g/\partial z = 0$ and 
$h'_1\partial g/\partial x + 
h'_2\partial g/\partial y + 
h'_3\partial g/\partial z = 0$.  
Let $D$ be the divisor of $S$ defined by 
$\partial g/\partial x = 0$. 
Then the above 2-form can be thought of 
as an element of 
$H^0(S^0, {\cal F}\vert_{S^0}(D))$. 
Now let us consider $\omega_2$ as an 
element of $H^0(S^0, {\cal F}\vert_{S^0}(D))$ 
by the injection 
$H^0(S^0, {\cal F}\vert_{S^0}) \to 
H^0(S^0, {\cal F}\vert_{S^0}(D))$. We put 

$$ \omega_3 := \omega_2 - 
\Sigma_{1 \leq j \leq n-3}(h^{(j)}_2 
dz/(\partial g/\partial x) - h^{(j)}_3 
dy/(\partial g/\partial x)) \wedge 
ds_j  -  \{h'_2 dz/(\partial g/\partial x) - 
h'_3 dy/(\partial g/\partial x\} \wedge dt. $$ 

By the exact sequence 

$$  0 \to H^0(S^0, 
(p^0)^*\Omega^2_{\Delta^{n-2}}\vert_{S^0}(D)) 
\to H^0(S^0, F\vert_{S^0}(D)) \to 
H^0(S^0, \Omega^1_{S^0}(D)\otimes 
(p^0)^*\Omega^1_{\Delta^{n-2}}) $$

we have $\omega_3 \in H^0(S^0, 
(p^0)^*\Omega^2_{\Delta^{n-2}}\vert_{S^0}(D))$. 
We write 

$$  \omega_3 = \Sigma_{1 \leq j \leq n-3}k^{(j)} 
ds_j \wedge dt/(\partial g/\partial x) + 
\Sigma_{1 \leq l, m \leq n-3}k^{(l,m)} 
ds_l \wedge ds_m/\partial g/\partial x) $$ 

with holomorphic functions $h^{(j)}$, 
$h^{(l,m)}$ on $S^0$. We shall check the condition 
that $\omega_2 \in 
H^0(S^0, {\cal F}\vert_{S^0})$. 
In particular, $\omega_2 \in 
\Gamma (S^0 \setminus \{\partial g/\partial y = 0 \}, 
{\cal F}\vert_{S^0})$ and $\omega_2 \in 
\Gamma (S^0 \setminus \{\partial g/\partial z = 0 \}, 
{\cal F}\vert_{S^0})$. 
In order to check these conditions, 
we substitute $dy$ by 

$$  -1/(\partial g/\partial y)\{\partial g/\partial x dx 
+ \partial g/\partial z dz + \partial g/\partial t dt\}  $$  

or substitute $dz$ by 

$$  -1/(\partial g/\partial z)\{\partial g/\partial x dx 
+ \partial g/\partial y dy + 
\partial g/\partial t dt\}.  $$

Note that the terms of $ds_j$ do not 
appear because 
$(\partial g/\partial s_i)\vert_{t=0} 
\equiv 0$ by the property (a).  
                   
By using the equations $h^{(j)}_1\partial g/\partial x 
+ h^{(j)}_2\partial g/\partial y + 
h^{(j)}_3\partial g/\partial z = 0$ and 
$h'_1\partial g/\partial x + 
h'_2\partial g/\partial y + 
h'_3\partial g/\partial z = 0$, 
we finally have the following 
two expressions of $\omega_2$: \vspace{0.4cm}

{\bf (y)}: $  \omega_2 = 
\Sigma_{1 \leq j \leq n-3} 
\{h^{(j)}_3dx/(\partial g/\partial y) - 
h^{(j)}_1dz/(\partial g/\partial y)\} 
\wedge ds_j $ \vspace{0.2cm}

$ + \{h'_3dx/(\partial g/\partial y) - 
h'_1dz/(\partial g/\partial y)\} 
\wedge dt $ \vspace{0.2cm}

$ + \Sigma_{1 \leq j \leq n-3} 
\{k^{(j)}/(\partial g/\partial x) - 
h^{(j)}_3(\partial g/\partial t)/(\partial g/\partial x)
(\partial g/\partial y)\}ds_j \wedge dt $ \vspace{0.2cm}

$ + \Sigma_{1 \leq l, m \leq n-3}k^{(l,m)} 
ds_l \wedge ds_m/(\partial g/\partial x)$.  
\vspace{0.4cm}

{\bf (z)}: $  \omega_2 = 
\Sigma_{1 \leq j \leq n-3} 
\{-h^{(j)}_2dx/(\partial g/\partial z) + 
h^{(j)}_1dy/(\partial g/\partial z)\} 
\wedge ds_j $ \vspace{0.2cm}

$ + \{-h'_2dx/(\partial g/\partial z) + 
h'_1dy/(\partial g/\partial z)\} 
\wedge dt $ \vspace{0.2cm}

$ + \Sigma_{1 \leq j \leq n-3} 
\{k^{(j)}/(\partial g/\partial x) + 
h^{(j)}_2(\partial g/\partial t)/(\partial g/\partial x)
(\partial g/\partial z)\}ds_j \wedge dt $ \vspace{0.2cm}

$ + \Sigma_{1 \leq l, m \leq n-3}k^{(l,m)} 
ds_l \wedge ds_m/(\partial g/\partial x)$.  \vspace{0.4cm}
 
{\bf Lemma (2.3)}. $k^{(j)} \in (x)+m^2{\cal O}_{S,0}$ 
$(1 \leq j \leq n-3)$, and 
$k^{(l,m)} \in x{\cal O}_{S,0}$ 
$(1 \leq l, m \leq n-3)$.  \vspace{0.12cm}

{\em Proof}. $k^{(l,m)}/(\partial g/\partial x)$ 
should be regular on $S^0$. 
This implies that $k^{(l,m)} \in x{\cal O}_{S,0}$ 
because $\partial g/\partial x = 2x$. 
We shall prove that $k^{(j)} \in m{\cal O}_{S,0}$ 
by using property (b) of $g(x,y,z,s,t)$ 
in each case $A_n$, $D_n$, $E_6$, $E_7$ or 
$E_8$. \vspace{0.15cm}

$(A_n)$: Since 
$\partial g/\partial t = 0$, $k^{(j)} 
\in x{\cal O}_{S,0}$ by the same reason as 
$h^{(l,m)}$ are so. \vspace{0.12cm}

$(D_n)$: $\omega\vert_{S^0}$ should be 
regular on $S \setminus \{\partial g/\partial z = 0\}$. 
In the expression {\bf (z)}, 
we denote by $G^{(j)}$ the coefficients of 
$ds_j \wedge dt$. $G^{(j)}$ is regular on 
$S \setminus \{\partial g/\partial z = 0\}$. 
In our case, $\partial g/\partial x = 2x$, 
$\partial g/\partial z = y^2 + (n-1)z^{n-2}$ 
and $\partial g/\partial t = cz^{n-2}$ 
(c: constant). Since it is easily checked 
that $\{\partial g/\partial x = 0 \}$ 
and $\{\partial g/\partial z = 0 \}$ 
have no common irreducible components 
on $S$, this implies that 
$k^{(j)}(\partial g/\partial z) + 
h^{(j)}_2(\partial g/\partial t) 
\in (\partial g/\partial x){\cal O}_{S,0}$. 

More explicitly, we have 
$\{y^2 + (n-1)z^{n-2}\}k^{(j)} + 
cz^{n-2}h^{(j)}_2 \in x{\cal O}_{S,0}$. 
We see at first that $k^{(j)} 
\in m{\cal O}_{S,0}$ because $y^2$ 
does not appear in the quadratic part 
of $g(x,y,z,0,0) = x^2 + y^2z +z^{n-1}$. 

Next the linear part of $k^{(j)}$ does 
not contain $y$ because any element of 
the form $\alpha y^3$ + (other terms), 
with $\alpha \in {\bold C}^*$, is not 
contained in the ideal 
$(x, z^{n-2}, x^2 + y^2z + z^{n-1})$ of 
${\bold C}\{x,y,z\}$. 

Finally the linear part of 
$k^{(j)}$ does not contain $z$ because $z$ 
does not appear in the linear part of 
$h^{(j)}_2$ (by an explicit calculation 
using the equation $h^{(j)}_1\partial g/\partial x 
+ h^{(j)}_2\partial g/\partial y + 
h^{(j)}_3\partial g/\partial z = 0$ in 
${\cal O}_{S,0}$) and therefore any 
element of the form $\alpha 
\{y^2z + (n-1)z^{n-1}\}$ + (other terms), 
with $\alpha \in {\bold C}^*$, is not 
contained in the ideal 
$(x^2 + y^2z + z^{n-1}, x, h^{(j)}_2z^{n-2})$ 
of ${\bold C}\{x,y,z\}$. \vspace{0.12cm}    
  
$(E_6)$: $\omega\vert_{S^0}$ should be 
regular on $S \setminus \{\partial g/\partial y = 0\}$. 
In the expression {\bf (y)}, 
we denote by $F^{(j)}$ the coefficients of 
$ds_j \wedge dt$. $F^{(j)}$ is regular on 
$S \setminus \{\partial g/\partial y = 0\}$. 
In our case, $\partial g/\partial x = 2x$, 
$\partial g/\partial y = 3y^2 $ and 
$\partial g/\partial t = cyz^2$ (c: constant). 
Since it is easily checked that 
$\{\partial g/\partial x = 0 \}$ and 
$\{\partial g/\partial y = 0 \}$ 
have no common irreducible components on $S$, 
this implies that $k^{(j)}(\partial g/\partial y) 
- h^{(j)}_3(\partial g/\partial t) 
\in (\partial g/\partial x){\cal O}_{S,0}$. 
More explicitly, we have $3yk^{(j)} - 
cz^2h^{(j)}_3 \in x{\cal O}_{S,0}$. 
We have $k^{(j)} \in m{\cal O}_{S,0}$ 
because $g(x,y,z,0,0)$ has no linear terms. 
Moreover, $y$ and $z$ do not appear 
in the linear part of $k^{(j)}$ because 
any element of the form $\alpha y^2$ + 
(other terms) or of the form $\alpha yz$ + 
(other terms), $(\alpha \in {\bold C}^*)$, 
is not contained in the ideal 
$(x, x^2 + y^3 + z^4, z^2)$ of 
${\bold C}\{x,y,z\}$. \vspace{0.12cm}

$(E_7)$: $\omega\vert_{S^0}$ should be regular 
on $S \setminus \{\partial g/\partial y = 0\}$. 
In the expression {\bf (y)}, we denote by 
$F^{(j)}$ the coefficients of $ds_j 
\wedge dt$. $F^{(j)}$ is regular on 
$S \setminus \{\partial g/\partial y = 0\}$. 
In our case, $\partial g/\partial x = 2x$, 
$\partial g/\partial y = 3y^2 + z^3$ and 
$\partial g/\partial t = cz^4$ (c: constant). 
Since it is easily checked that 
$\{\partial g/\partial x = 0 \}$ and 
$\{\partial g/\partial y = 0 \}$ 
have no common irreducible components on 
$S$, this implies that 
$k^{(j)}(\partial g/\partial y) - 
h^{(j)}_3(\partial g/\partial t) 
\in (\partial g/\partial x){\cal O}_{S,0}$. 
More explicitly, we have 
$(3y^2 + z^3)k^{(j)} - cz^4h^{(j)}_3 
\in x{\cal O}_{S,0}$. We have $k^{(j)} 
\in m{\cal O}_{S,0}$ because the term 
$y^2$ is not contained in the 
quadratic part of $g(x,y,z,0,0)$. 
Moreover, $y$ and $z$ do not appear in 
the linear part of $k^{(j)}$ 
because any element of the form 
$\alpha (3y^3 + yz^3)$ + (other terms) 
or of the form  $\alpha y^2z$ + 
(other terms), $(\alpha \in {\bold C}^*)$, 
is not contained in the ideal 
$(x, h^{(j)}_3z^4, x^2 + y^3 + yz^3)$ 
of ${\bold C}\{x,y,z\}$. \vspace{0.12cm}
    
$(E_8)$: $\omega\vert_{S^0}$ should be 
regular on $S \setminus 
\{\partial g/\partial y = 0\}$. 
In the expression {\bf (y)}, 
we denote by $F^{(j)}$ the coefficients of 
$ds_j \wedge dt$. $F^{(j)}$ is regular on 
$S \setminus \{\partial g/\partial y = 0\}$. 
In our case, $\partial g/\partial x = 2x$, 
$\partial g/\partial y = 3y^2 $ and 
$\partial g/\partial t = cyz^3$ (c: constant). 
Since it is easily checked that 
$\{\partial g/\partial x = 0 \}$ and 
$\{\partial g/\partial y = 0 \}$ 
have no common irreducible components on 
$S$, this implies that 
$k^{(j)}(\partial g/\partial y) - 
h^{(j)}_3(\partial g/\partial t) 
\in (\partial g/\partial x){\cal O}_{S,0}$. 
More explicitly, we have $3yk^{(j)} - 
cz^3h^{(j)}_3 \in x{\cal O}_{S,0}$. 
We have $k^{(j)} \in m{\cal O}_{S,0}$ 
because $g(x,y,z,0,0)$ has no linear terms. 
Moreover, $y$ and $z$ do not appear in 
the linear part of $k^{(j)}$ because 
any element of the form $\alpha y^2$ + 
(other terms) or of the form $\alpha yz$ + 
(other terms), $(\alpha \in {\bold C}^*)$, 
is not contaied in the ideal 
$(x, z^3, x^2 + y^3 + z^5)$ of 
${\bold C}\{x,y,z\}$. \vspace{0.15cm}

{\bf Step 3.} We shall write $\omega\vert_{S^0}$ 
explicitly and calculate 
$\wedge^{n/2}\omega\vert_{S^0}$. 
Summing up $\omega_1$ and $\omega_2$, 
we have \vspace{0.3cm}

$ \omega\vert_{S^0} = \{\overline{f}_2 + h'_2\}dz 
\wedge dt/(\partial g/\partial x) - 
\{\overline{f}_3 + h'_3\}dy \wedge 
dt/(\partial g/\partial x)$ \vspace{0.2cm}

$ + \overline{f}_4dy \wedge 
dz/(\partial g/\partial x)$ \vspace{0.2cm}

$ + \Sigma_{1 \leq j \leq n-3}\{h^{(j)}_2 
dz/(\partial g/\partial x) - h^{(j)}_3 
dy/(\partial g/\partial x)\} \wedge ds_j $ 
\vspace{0.2cm}

$ + \Sigma_{1 \leq j \leq n-3}k^{(j)} 
ds_j \wedge dt/(\partial g/\partial x) + 
\Sigma_{1 \leq l, m \leq n-3}{k'}^{(l,m)} 
ds_l \wedge ds_m. $\vspace{0.2cm}

By Lemmas (2.1), (2.2) and (2.3) the functions 
$\overline{f}_i$, $h^{(j)}_i$, $h'_i$ are 
all contained in $m{\cal O}_{S,0}$, and 
$k^{(j)} \in (x)+m^2{\cal O}_{S,0}$. 
By Lemma (2.3), ${k'}^{(l,m)} := 
k^{(l,m)}/(\partial g/\partial x)$ 
is regular. 

An explicit calculation shows that

$$ \wedge^{n/2}\omega\vert_{S^0} = 
B\cdot dy \wedge dz \wedge dt \wedge ds_1
\wedge ... \wedge ds_{n-3} $$

where $B$ is of the following form 
\vspace{0.2cm}

$ (\overline f_2 + 
h'_2)/(\partial g/\partial x)\Sigma_{1 \le j \le n-3}
\{h^{(j)}_3/(\partial g/\partial x)
\cdot$(regular function)\} \vspace{0.2cm}

$+ (\overline f_3 + h'_3)/(\partial g/\partial
x)\Sigma_{1 \le j \le n-3}
\{h^{(j)}_2/(\partial g/\partial x)
\cdot$(regular function)\}
\vspace{0.2cm}

$+ \overline f_4/(\partial g/\partial x)
\Sigma_{1 \le j \le n-3}\{k^{(j)}/
(\partial g/\partial x)
\cdot$(regular function)\} \vspace{0.2cm}

$+ \Sigma_{1 \leq p,q,r \leq n-3}h^{(p)}_2/
(\partial g/\partial x)\cdot
h^{(q)}_3/(\partial g/\partial x)
\cdot k^{(r)}/(\partial g/\partial
x)\cdot$(regular function). \vspace{0.3cm}

When $S$ is of type $A_n$, then 
$k^{(p)} \in x{\cal O}_{S,0}$ by the proof of
Lemma (2.3). Hence $B = C/x^2$ for a suitable 
$C \in m^2$, where $m$ is the
maximal ideal of ${\cal O}_{S,0}$. 
If $\wedge^{n/2}\omega\vert_{S^0}$ is
nowhere-vanishing (i.e. 
$\wedge^{n/2}\omega\vert_{S^0}$ is a generator of
$\omega_{\cal H}\vert_S$), 
then we should have $B = (unit)/x$. This means
that $C = (unit)x$, hence $x \in m^2$, 
which is absurd.

When $S$ is one of other types. 
Then $B = C/x^3$ for a suitable $C \in xm^2
+ m^4$. If $\omega\vert_{S^0}$ is 
nowhere-vanishing, then we should have $B
= (unit)/x$. This means that 
$C = (unit)x^2$, hence $x^2 \in xm^2 + m^4$.
However this is absurd because $g(x,y,z,0,0)$ 
always contains a cubic term
which cannot be divided by $x$. 
\vspace{0.2cm}

\begin{center}
Department of Mathematics, Graduate school 
of Sience, Osaka University, Toyonaka 560, Japan
\end{center}

\end{document}